\documentclass[a4paper,11pt]{article}
\usepackage{amsfonts}
\usepackage{amstext}
\usepackage{amsthm}
\usepackage{amsmath}
\usepackage{amssymb}
\usepackage{latexsym}
\usepackage{graphicx}
\usepackage[latin1]{inputenc}
\newcommand{\R}{\Bbb{R}}

\newtheorem{theorem}{Theorem}[section]

\newcommand{\n}{\noindent}

\begin{document}

\title{Maximum and comparison principles for degenerate elliptic systems and some applications
}

\author{\textbf{Edir Junior Ferreira Leite \footnote{\textit{E-mail addresses}:
edirjrleite@ufv.br (E.J.F. Leite)}}\\ {\small\it Departamento de Matem\'{a}tica,
Universidade Federal de Vi\c{c}osa,}\\ {\small\it CCE, 36570-900, Vi\c{c}osa, MG, Brazil}}

\date{}{%{\it Preprint - December 12, 2009}}

\maketitle

\markboth{abstract}{abstract}
\addcontentsline{toc}{chapter}{abstract}

\hrule \vspace{0,2cm}

\n {\bf Abstract}

In this paper we develop a detailed study on maximum and comparison principles related to the following nonlinear eigenvalue problem
\[
\left\{
\begin{array}{llll}
-\Delta_p u = \lambda a(x)\vert v\vert^{\beta_1-1}v & {\rm in} \ \ \Omega;\\
-\Delta_q v = \mu b(x)\vert u\vert^{\beta_2-1}u & {\rm in} \ \ \Omega;\\
u= v=0 & {\rm on} \ \ \partial\Omega,
\end{array}
\right.
\]
where $p,q\in (1,\infty)$, $\beta_1, \beta_2 > 0$ satisfy $\beta_1\beta_2=(p-1)(q-1)$, $\Omega\subset\R^n$ is a bounded domain with $C^2$-boundary, $a,b\in L^\infty(\Omega)$ are given functions, both assumed to be strictly positive on compact subsets of $\Omega$, and $\Delta_p$ and $\Delta_q$ are quasilinear elliptic operators, stand for $p$-Laplacian and $q$-Laplacian, respectively. We classify all couples $(\lambda, \mu) \in \R^2$ such that both the (weak and strong) maximum and comparison principles corresponding to the above system hold in $\Omega$. Explicit lower bounds for principal eigenvalues of this system in terms of the measure of $\Omega$ are also proved. As application, given $\lambda,\mu\geq 0$ we measure explicitly how small has to be $\vert \Omega\vert$ so that weak and strong maximum principles associated to the above problem hold in $\Omega$.

\vspace{0.5cm}
\hrule\vspace{0.2cm}

\section{Introduction and statements}

In this paper we study maximum and comparison principles related to the following system:
\begin{equation}\label{1.3}
\left\{
\begin{array}{llll}
-\Delta_pu = \lambda a(x)\vert v\vert^{\beta_1-1}v & {\rm in} \ \ \Omega;\\
-\Delta_qv = \mu b(x)\vert u\vert^{\beta_2-1}u & {\rm in} \ \ \Omega;\\
u= v=0 & {\rm on} \ \ \partial\Omega,
\end{array}
\right.
\end{equation}
where $\Omega\subset\R^n$ is a bounded domain with $C^2$-boundary (not necessarily connected), $p,q\in (1,\infty)$, $a,b\in L^\infty(\Omega)$ are given functions satisfying 
\[
\underset{x\in\Omega}{\mathrm{ess }\inf}\ a(x)>0\ \ \text{ and }\ \ \underset{x\in\Omega}{\mathrm{ess }\inf}\ b(x)>0,
\]
$\beta_1, \beta_2 > 0$ with $\beta_1 \beta_2 = (p-1)(q-1)$, and $(\lambda, \mu) \in \R^2$. The $p$-Laplacian is defined by
\[
\Delta_pu:=\operatorname{div}\, (\vert\nabla u\vert^{p-2}\nabla u)
\] 
for any $u\in W^{1,p}_0(\Omega)$ with values $\Delta_pu\in W^{-1,\frac{p}{p-1}}(\Omega)$, the dual space of $W^{1,p}_0(\Omega)$.

Existence, nonexistence and uniqueness of nontrivial solutions to the system (\ref{1.3}) have been widely investigated during the three last decades for $p = q = 2$ and, more generally, for $p,q\in (1,\infty)$. For $p = q = 2$, we refer for instance to \cite{ClDfMi}, \cite{DfFe}, \cite{FeMa}, \cite{HuvdV}, \cite{Mi1} and \cite{SZ1}, where in particular notions of sub-superlinearity, sub-supercriticality and criticality  have been introduced. Still in the first part, the eigenvalue problem, i.e., $\beta_1\beta_2 = 1$, was completely studied in \cite{Mo4}. For $p,q\in (1,\infty)$, we refer to \cite{2} when $\beta_1\beta_2>(p-1)(q-1)$ and \cite{Cuesta}  when $\beta_1\beta_2=(p-1)(q-1)$.

The connection between principal eigenvalues and maximum principles have been investigated in \cite{AFM, Am2, AnLo, CauLo, LoMo} for cooperative systems and in \cite{Sw} for non-cooperative systems (see also \cite{Lo1} for a more complete discussion) and more recently in \cite{LM}, where system (\ref{1.3}) was analyzed in the special case when $p = q = 2$, however, instead of $\Delta$, a general second order elliptic operator was considered. For the nonlocal context, we refer to \cite{EM1}.

Here we extend the results of \cite{LM} for $p,q>1$, that is, we establish the connection between principal spectral curves for systems (\ref{1.3}) and maximum and comparison principles related. For this, we shall present a bit of notation. Note that, given any $f\in L^\infty(\Omega)$, there exists a unique weak solution $u\in W^{1,p}_0(\Omega)$ of the classical problem

\begin{equation}\label{P1}
\left\{
\begin{array}{rrll}
-\Delta_p u &=& f(x) & {\rm in} \ \ \Omega;\\
u &=& 0  & {\rm on} \ \ \partial \Omega.
\end{array}
\right.
\end{equation}
Notice that, $u\in C^{1,\alpha}(\overline{\Omega})$ for some $\alpha\in (0,1)$ (see \cite{7, GS, 15, 24}). We denote $X:=\left[C^1_0(\overline{\Omega})\right]^2$, $X_+:=\{(u,v)\in X:u\geq 0\text{ and }v\geq 0\text{ in }\Omega\}$, and $\stackrel{\circ}{X}_+$ is the topological interior of $X_+$ in $X$. Note that, $\stackrel{\circ}{X}_+$ is nonempty and characterized by $(u,v)\in\ \stackrel{\circ}{X}_+$ if, and only if, $(u,v)\in X$ satisfies:
\[
u,v>0\ \text{in}\ \Omega\ \ \ \text{and}\ \ \ \frac{\partial u}{\partial \nu},\frac{\partial v}{\partial \nu}<0\ \text{on}\ \partial\Omega,
\]
where $\nu\equiv \nu(x_0)$ denotes the exterior unit normal to $\partial\Omega$ at $x_0\in\partial\Omega$ (see \cite{Cuesta}).

As is well known, the operator $\Delta_p$ satisfies the weak maximum principle, that is, for any weak solution $u \in W^{1,p}_0(\Omega)$ to 
\[
\left\{
\begin{array}{rrll}
-\Delta_p u &=& f(x) & {\rm in} \ \ \Omega;\\
u &\geq& 0  & {\rm on} \ \ \partial \Omega,
\end{array}
\right.
\]
with $f\in L^{\infty}(\Omega)$ and $f\geq 0$ in $\Omega$, one has $u \geq 0$ in $\Omega$ and $\Delta_p$ also satisfies the strong maximum principle, i.e., moreover $u > 0$ in $\Omega$ whenever $f\not\equiv 0$ in $\Omega$ (see \cite{GS}, \cite{23} and \cite{25}).

Let $(u,v)$ in $W^{1,p}_0(\Omega)\times W^{1,q}_0(\Omega)$. The weak formulation of the system (\ref{1.3}) is given by
\[
\lambda\int_{\Omega}a(x)\vert v\vert^{\beta_1-1}vw dx= \int_{\Omega}\vert\nabla u\vert^{p-2}\nabla u\nabla w dx,\ \ \ \ \forall\ w \in C^1_0(\Omega)
\]
and
\[
\mu\int_{\Omega}b(x)\vert u\vert^{\beta_2-1}uw dx= \int_{\Omega}\vert\nabla v\vert^{q-2}\nabla v\nabla w dx,\ \ \ \ \forall\ w \in C^1_0(\Omega).
\]

A couple $(\lambda, \mu) \in \R_+^*\times\R_+^*=(0,\infty)^2$ is said to be an eigenvalue of \eqref{1.3} if the system admits a nontrivial weak solution $(\varphi,\psi)$ in $W^{1,p}_0(\Omega)\times W^{1,q}_0(\Omega)$ which is called an eigenfunction associated to $(\lambda, \mu)$. We say that $(\lambda, \mu)$ is a principal eigenvalue if admits a positive eigenfunction $(\varphi,\psi)$; i.e., $\varphi$ and $\psi$ are positive in $\Omega$. We also say that $(\lambda, \mu)$ is simple in $\stackrel{\circ}{X}_+$ if for any eigenfunctions $(\varphi,\psi),(\tilde{\varphi}, \tilde{\psi})\in\ \stackrel{\circ}{X}_+$, there exists $\rho>0$ such that $\tilde{\varphi} = \rho \varphi$ and $\tilde{\psi} = \rho\mu^{\frac{1}{\beta_2}} \psi$ in $\Omega$.

The existence of principal eigenvalues of \eqref{1.3} and some of their qualitative properties were treated in Cuesta and Takác \cite{Cuesta}. Namely, they proved that the set formed by these principal eigenvalues is given by the following smooth curve
\[
\mathcal{C}_1:=\left\{(\lambda,\mu)\in (\R_+^*)^2:\lambda^{\frac{1}{\sqrt{\beta_1(p-1)}}}\mu^{\frac{1}{\sqrt{\beta_2(q-1)}}}=\Lambda' \right\},
\]
for some $\Lambda'>0$, which satisfies:

\begin{itemize}

\item[(i)] $(\lambda, \mu) \in \R_+\times\R_+$ is a principal eigenvalue of the system (\ref{1.3}) if, and only if, $(\lambda,\mu)\in \mathcal{C}_1$;

\item[(ii)] The curve $\mathcal{C}_1$ is simple in $\stackrel{\circ}{X}_+$, that is, $(\lambda,\mu)$ is simple in $\stackrel{\circ}{X}_+$ for all $(\lambda,\mu)\in \mathcal{C}_1$;

\item[(iii)] Let $(\varphi,\psi)\in X$ be an eigenfunction corresponding to $(\lambda,\mu)\in \mathcal{C}_1$. Therefore either  $(\varphi,\psi)\in\ \stackrel{\circ}{X}_+$ or $(-\varphi,-\psi)\in\ \stackrel{\circ}{X}_+$.
\end{itemize}

Recently, in \cite{L} the author established another properties satisfied by principal curve $\mathcal{C}_1$ such as local isolation, monotonicity on $\Omega$ and monotonicity and continuity of the principal eigenvalue with respect to the weight functions $a$ and $b$.

By weak maximum principle, denoted by {\bf (WMP)}, we mean that for any weak solution $(u,v)\in W^{1,p}_0(\Omega)\times W^{1,q}_0(\Omega)$ of the system 
\begin{equation}\label{supersol}
\left\{
\begin{array}{llll}
-\Delta_p u = \lambda a(x)\vert v\vert^{\beta_1-1}v +f(x) & {\rm in} \ \ \Omega;\\
-\Delta_qv = \mu b(x)\vert u\vert^{\beta_2-1}u + g(x) & {\rm in} \ \ \Omega;\\
u= v=0 & {\rm on} \ \ \partial\Omega,
\end{array}
\right.
\end{equation}
with $f,g\in L^{\infty}(\Omega)$ and $f,g\geq 0$ in $\Omega$, verifies $u, v \geq  0$ in $\Omega$. Besides, if at least, $u$ or $v$ is positive in $\Omega$ whenever $f+g\not\equiv 0$ in $\Omega$, we say that the strong maximum principle, denoted by {\bf (SMP)}, corresponding to \eqref{1.3} holds in $\Omega$. When $\lambda, \mu > 0$, {\bf (SMP)} can be rephrased as $u, v > 0$ in $\Omega$ whenever $f+g\not\equiv 0$ in $\Omega$. In the case that $u > 0$ ($v > 0$) in $\Omega$, we get $\frac{\partial}{\partial \nu}u(x_0) < 0$ $\left(\frac{\partial}{\partial \nu}v(x_0) < 0\right)$ for all $x_0 \in \partial \Omega$.

We are ready to classify completely in terms of the principal curve $\mathcal{C}_1$ the set of couples $(\lambda, \mu) \in \R^2$ such that {\bf (WMP)} and {\bf (SMP)} hold in $\Omega$.

Namely:

\begin{theorem} \label{MP}
Let $(\lambda, \mu) \in \R^2$ and $\mathcal{R}_1$ be the open region in the first quadrant below $\mathcal{C}_1$. The following assertions are equivalent:

\begin{itemize}
\item[{\rm (i)}] $(\lambda, \mu) \in \overline{\mathcal{R}_1} \setminus \mathcal{C}_1$;
\item[{\rm (ii)}] {\bf (WMP)} corresponding to \eqref{1.3} holds in $\Omega$;
\item[{\rm (iii)}] {\bf (SMP)} corresponding to \eqref{1.3} holds in $\Omega$.
\end{itemize}
\end{theorem}

Note that the sets 
\[
\mathcal{R}_1=\left\{(\lambda,\mu)\in (\R_+^*)^2:\lambda^{\frac{1}{\sqrt{\beta_1(p-1)}}}\mu^{\frac{1}{\sqrt{\beta_2(q-1)}}}<\Lambda' \right\}
\]
and $\overline{\mathcal{R}_1} \setminus \mathcal{C}_1 = \mathcal{R}_1 \cup \{ (\lambda,0):\ \lambda \geq 0\} \cup \{ (0,\mu):\ \mu \geq 0\}$ are unbounded and are depicted in Figure 1.

\begin{figure}[ht]
\centering
\includegraphics[scale=1]{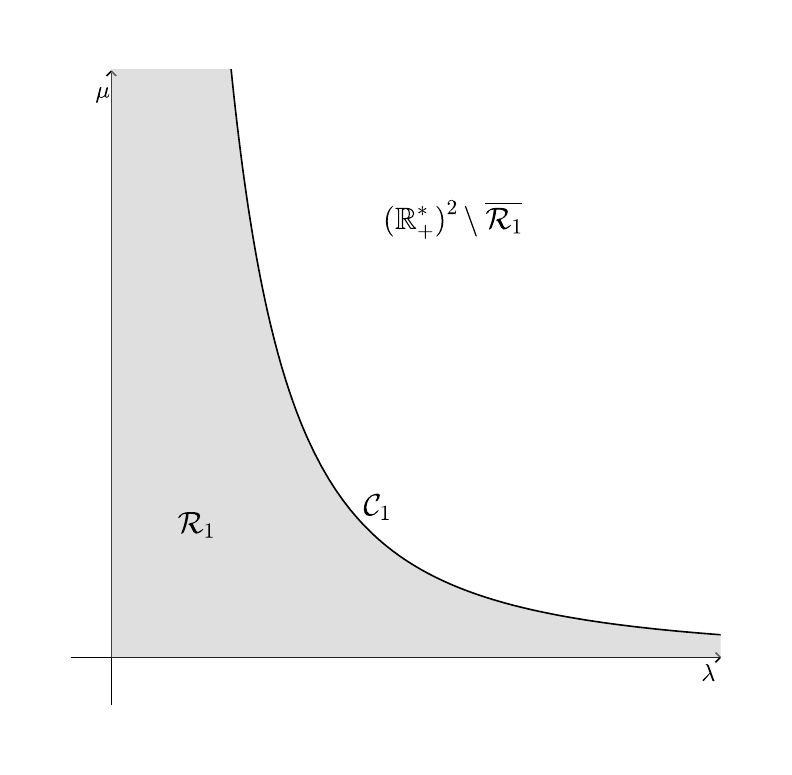}
\caption{The principal curve $\mathcal{C}_1$.}
\end{figure}

We notice that weak and strong comparison principles, denoted respectively by {\bf (WCP)} and {\bf (SCP)}, are very important tools for establish the uniqueness and positivity of solutions for elliptic problems, among others, to certain counterparts of (\ref{1.3}). We say that {\bf (WCP)} holds in $\Omega$ if, for any weak solutions $(u,v)$ and $(z,w)$ in $W^{1,p}_0(\Omega)\times W^{1,q}_0(\Omega)$ to the following systems, respectively:
\begin{equation}\label{comp1}
\left\{
\begin{array}{llll}
-\Delta_p u = \lambda a(x)\vert v\vert^{\beta_1-1}v +f_1(x) & {\rm in} \ \ \Omega;\\
-\Delta_qv = \mu b(x)\vert u\vert^{\beta_2-1}u + g_1(x) & {\rm in} \ \ \Omega;\\
u= v=0 & {\rm on} \ \ \partial\Omega,
\end{array}
\right.
\end{equation}
\begin{equation}\label{comp2}
\left\{
\begin{array}{llll}
-\Delta_p z = \lambda a(x)\vert w\vert^{\beta_1-1}w +f_2(x) & {\rm in} \ \ \Omega;\\
-\Delta_qw = \mu b(x)\vert z\vert^{\beta_2-1}z + g_2(x) & {\rm in} \ \ \Omega;\\
z= w=0 & {\rm on} \ \ \partial\Omega,
\end{array}
\right.
\end{equation}
with $f_1,f_2,g_1,g_2\in L^\infty(\Omega)$ and $0\leq f_1 \leq f_2$ and $0\leq g_1 \leq  g_2$ in $\Omega$, one has $u \leq z$ and $v \leq w$ in $\Omega$. If, in addition, at least, $u < z$ or $v < w$ in $\Omega$ whenever $f_1+g_1\not\equiv f_2+g_2$ in $\Omega$, we say that {\bf (SCP)} corresponding to (\ref{1.3}) holds in $\Omega$. When $\lambda, \mu >0$, {\bf (SCP)} in $\Omega$ can be rewritten as $u < z$ and $v < w$ in $\Omega$  whenever $f_1+g_1\not\equiv f_2+g_2$ in $\Omega$. In the case that $u < z$ ($v < w$) in $\Omega$, we clearly have $\frac{\partial}{\partial \nu}u(x_0) > \frac{\partial}{\partial \nu}z(x_0)$ $\left(\frac{\partial}{\partial \nu}v(x_0) > \frac{\partial}{\partial \nu}w(x_0)\right)$ for all $x_0 \in \partial \Omega$.

Our next theorem characterizes completely the {\bf (WCP)} and {\bf (SCP)} corresponding to (\ref{1.3}) in terms of the smooth curve $\mathcal{C}_1$. Precisely:

\begin{theorem} \label{CP} Let $\mathcal{R}_1$ be as in Theorem \ref{MP} and $(\lambda, \mu) \in \R^2$. The following assertions are equivalent:

\begin{itemize}
\item[{\rm (i)}] $(\lambda, \mu) \in \overline{\mathcal{R}_1} \setminus \mathcal{C}_1$;
\item[{\rm (ii)}] {\bf (WCP)} corresponding to \eqref{1.3} holds in $\Omega$;
\item[{\rm (iii)}] {\bf (SCP)} corresponding to \eqref{1.3} holds in $\Omega$.
\end{itemize}
\end{theorem}

The validity of weak and strong comparison principles for problems involving the $p$-Laplacian operator is usually very delicate (see \cite{CT, Cuesta, Da, DaSc, DaSc1, FGT, FHTT, FT, GS, GV, PS, 23}). For example,  the operator $\Delta_p$ satisfies the strong comparison principle only under the additional condition zero Dirichlet boundary values on $\partial\Omega$. So, we need assume the same condition on the {\bf (WMP)}, {\bf (SMP)}, {\bf (WCP)} and {\bf (SCP)} associated to the problem (\ref{1.3}). Moreover, we prove the {\bf (WCP)} and {\bf (SCP)} under the following additional condition: $ f_1\geq 0$ and $ g_1\geq 0$ in $\Omega$, when $(\lambda, \mu) \in \mathcal{R}_1$. 

Now, we characterize when such {\bf (WMP)} (or {\bf (SMP)}) corresponding to \eqref{1.3} is satisfied in domains $\Omega$ of small Lebesgue measure. For this, we shall obtain an explicit lower estimate of $\Lambda'$ in terms of the Lebesque measure of $\Omega$.

In Theorem 2.6 of \cite{BeNiVa}, Berestycki, Nirenberg and Varadhan established a lower estimate for principal eigenvalues corresponding to the problem, in the scalar context, involving linear second order elliptic operators. In  Theorem 5.1 of \cite{Lo}, López-Gómez obtained an explicit lower estimate for such principal eigenvalues. Later, in Theorem 10.1 of \cite{CaLo} Cano-Casanova and López-Gómez extended this result for mixed boundary conditions. The key tool for the proof of the Theorem 2.6 of \cite{BeNiVa} is an ABP estimate for second order uniformly elliptic operators while in proofs of the Theorems 5.1 of \cite{Lo} and 10.1 of \cite{CaLo} is used in a crucial way the celebrated Faber-Krahn inequality of Faber \cite{faber} and Krahn \cite{krahn}. Recently, in \cite{L} the author obtained an explicit lower bounds for principal eigenvalues of the system (\ref{1.3}) via Faber-Krahn inequality for the first eigenvalue of $-\Delta_p$. Here, we will use an ABP estimate related to $p$-Laplacian operator with explicit constant (see Theorem 3 of \cite{ACI}). Namely:

\begin{theorem} \label{lower} Let $a,b\in L^\infty(\Omega)\cap C(\Omega)$ and $\mathcal{C}_1$ be the principal curve associated to \eqref{1.3}. Then 

\begin{eqnarray}
&& \Lambda' \geq \frac{\vert\Omega\vert^{-\left(\frac{1}{ns}+\frac{1}{nr}\right)}}{c(n,p)^{\frac{\beta_2}{s}}c(n,q)^{\frac{\beta_1}{r}}d^{\frac{\beta_2}{s}+\frac{\beta_1}{r}} \|a\|^{\frac{1}{r}}_{L^\infty(\Omega)} \|b\|^{\frac{1}{s}}_{L^\infty(\Omega)}},  \label{lb1} 
\end{eqnarray}
where $d:=\operatorname{diam }(\Omega)$, $r:=\sqrt{\beta_1(p-1)}$, $s:=\sqrt{\beta_2(q-1)}$, $|\cdot |$ stands for the Lebesgue measure of $\R^n$, 
\[
c(n,p):=\left(n\min\{1,p-1\}\vert B_1\vert^{\frac{1}{n}}\right)^{-\frac{1}{p-1}}
\]
and $B_1$ is the unit ball of $\R^n$. In particular,
\[
\lim_{\vert\Omega\vert\downarrow 0}\Lambda'=+\infty.
\]
\end{theorem}

Note that we get an explicit lower estimate of $\Lambda'$ in terms of the Lebesque measure of $\Omega$, diameter of $\Omega$, explicit constants $c(n,p), c(n,q)$ and the weighted functions $a, b \in L^\infty(\Omega)\cap C(\Omega)$. As an interesting consequence of Theorems \ref{MP} and \ref{lower}, we obtain the following characterization of maximum principles:

\begin{theorem} \label{sm} Let $r,s,c(n,p)$ and $c(n,q)$ be as in Theorem \ref{lower} and $a,b\in L^\infty(\Omega) \cap C(\Omega)$. Define $d:=\operatorname{diam}(\Omega)$ and
\[
\eta := \frac{1}{\left[\lambda^{\frac{1}{r}}\mu^{\frac{1}{s}}\left(c(n,p)^{\frac{\beta_2}{s}}c(n,q)^{\frac{\beta_1}{r}}d^{\frac{\beta_2}{s}+\frac{\beta_1}{r}} \|a\|^{\frac{1}{r}}_{L^\infty(\Omega)} \|b\|^{\frac{1}{s}}_{L^\infty(\Omega)}\right)\right]^{\frac{nrs}{r+s}}}.
\]
The following assertions are equivalent:

\begin{itemize}
\item[{\rm (i)}] $\lambda \geq 0$ and $\mu \geq 0$;
\item[{\rm (ii)}] {\bf (WMP)} corresponding to \eqref{1.3} holds in $\Omega$ provided that $|\Omega| < \eta$;
\item[{\rm (iii)}] {\bf (SMP)} corresponding to \eqref{1.3} holds in $\Omega$ provided that $|\Omega| < \eta$.
\end{itemize}
\end{theorem}

Finally, as application of Theorems \ref{MP} and \ref{CP}, we present a characterization in terms of the principal curve $\mathcal{C}_1$ the set of couples $(\lambda,\mu)\in\mathbb{R}^2$ such that the system (\ref{supersol}) admits a unique weak nonnegative solution $(u,v)$ for any pair $(f,g) \in (L^\infty(\Omega))^2$ of nonnegative functions. Precisely, we have:

\begin{theorem} \label{FHA}
Let $\mathcal{R}_1$ be as in Theorem \ref{MP} and $(\lambda,\mu) \in \R^2$. Then, the couple $(\lambda,\mu) \in \overline{\mathcal{R}_1} \setminus \mathcal{C}_1$ if, and only if, the system (\ref{supersol}) admits a unique weak solution $(u,v) \in W^{1,p}_0(\Omega)\times W^{1,q}_0(\Omega)$ and this satisfies $u, v \geq 0$ in $\Omega$ for any pair $(f,g) \in (L^\infty(\Omega))^2$ of nonnegative functions. Moreover, $(u,v) \in (C^{1,\alpha}(\overline{\Omega}))^2$ for some $\alpha\in (0,1)$.  
\end{theorem}

The rest of paper is organized into five sections. In Section 2 we show Theorem \ref{MP} by mean of maximum and comparison principles and Hopf's lemma related to the $p$-Laplacian operator. In Section 3 we characterize weak and strong comparison principles associated to the system (\ref{1.3}) stated in Theorem \ref{CP} by using Theorem \ref{MP} as a key tool. In Section 4 we get an explicit lower estimate of $\Lambda'$ in terms of the Lebesque measure of $\Omega$ stated Theorem \ref{lower}. The characterization of maximum principle in domains $\Omega$ of small Lebesgue measure stated in Theorem \ref{sm} is established in Section 5. Finally, in Section 6 we prove Theorem \ref{FHA} by using Theorems \ref{MP} and \ref{CP}. 

\section{Proof of Theorem \ref{MP}}

In this section we prove Theorem \ref{MP}, which classify completely the couples $(\lambda,\mu)\in\R^2$ such that {\bf (WMP)} and {\bf (SMP)} hold in $\Omega$. Notice that it suffices to show only that (i) $\Leftrightarrow$ (ii). In this case, we clearly have (ii) $\Leftrightarrow$ (iii). In fact, it is obvious that {\bf (SMP)} in $\Omega$ implies {\bf (WMP)} in $\Omega$. Conversely, assume that {\bf (WMP)} holds in $\Omega$ and let $(u,v)$ be a weak solution of the system (\ref{supersol}) with $f,g\in L^{\infty}(\Omega)$, $f,g\geq 0$ in $\Omega$ and $f+g\not\equiv 0$ in $\Omega$. Thus, $u, v \geq 0$ in $\Omega$ and, by (i), we obtain $\lambda, \mu \geq 0$. Then, the conclusion of {\bf (SMP)} follows, since $\Delta_p$ satisfies the strong maximum principle (see \cite{GS, 23, 25}).

In order to proof that {\bf (WMP)} in $\Omega$ leads to $(\lambda, \mu) \in \overline{\mathcal{R}_1} \setminus \mathcal{C}_1$, assume instead that $(\lambda, \mu) \not\in \overline{\mathcal{R}_1} \setminus \mathcal{C}_1$. Let $(\lambda, \mu) \in \mathcal{C}_1$ and $(\tilde{\varphi},\tilde{\psi})$ be a positive eigenfunction associated to $(\lambda, \mu)$. Then, $(-\tilde{\varphi},-\tilde{\psi})$ is a negative eigenfunction associated to $(\lambda, \mu)$ and so {\bf (WMP)} fails in $\Omega$.

Assume now that $(\lambda, \mu) \in \R^2$ is a fixed couple outside of $\overline{\mathcal{R}_1}$. If $(\lambda, \mu) \in (\R_+^*)^2$, we get $\lambda > \lambda_1$ and $\mu > \mu_1$, where $(\lambda_1, \mu_1)$ is a principal eigenvalue of (\ref{1.3}) with $\frac{\mu}{\lambda}=\frac{\mu_1}{\lambda_1}$. Denote by $(\varphi,\psi)$ a positive eigenfunction corresponding to $(\lambda_1, \mu_1)$. Thus, we derive

\[
\left\{
\begin{array}{llll}
-\Delta_p (-\varphi) - \lambda a(x) \vert -\psi\vert^{\beta_1-1}(-\psi) &=& - \lambda_1 a(x) \psi^{\beta_1} + \lambda a(x) \psi^{\beta_1}\\
&=& (\lambda - \lambda_1)a(x) \psi^{\beta_1} \geq (\not\equiv)\ 0  \ \ & {\rm in}\ \Omega;\\
-\Delta_q (-\psi) - \mu b(x)\vert -\varphi\vert^{\beta_2-1}(-\varphi) &=& -\mu_1 b(x) \varphi^{\beta_2} + \mu b(x) \varphi^{\beta_2} \\
&=& (\mu - \mu_1)b(x) \varphi^{\beta_2} \geq (\not\equiv)\ 0 \ \ & {\rm in}\ \Omega
\end{array}
\right. 
\]
and $-\varphi = 0 =-\psi$ on $\partial\Omega$. Since, $-\varphi,-\psi<0$ in $\Omega$, {\bf (WMP)} doesn't hold in $\Omega$.

Now, suppose that $\lambda < 0$. Then, there exists $(\lambda_1, \mu_1)\in \mathcal{C}_1$ with $\lambda_1 > 0$ small enough (and so $\mu_1>0$ large enough) so that $\lambda < - \lambda_1$ and $\mu > - \mu_1$. Then, we get

\[
\left\{
\begin{array}{llll}
-\Delta_p (-\varphi) - \lambda a(x) \psi^{\beta_1} &=& -\lambda_1 a(x) \psi^{\beta_1} - \lambda a(x) \psi^{\beta_1} \\
&=& -(\lambda + \lambda_1)a(x) \psi^{\beta_1} \geq (\not\equiv)\ 0 \ \ & {\rm in}\ \Omega; \\
-\Delta_q \psi - \mu b(x)\vert -\varphi\vert^{\beta_2-1}(-\varphi) &=& \mu_1 b(x) \varphi^{\beta_2} + \mu b(x) \varphi^{\beta_2} \\
&=& (\mu+\mu_1) b(x) \varphi^{\beta_2} \geq (\not\equiv)\ 0\ \ & {\rm in}\ \Omega
\end{array}
\right. 
\]
and $-\varphi = 0 =\psi$ on $\partial\Omega$. However, $-\varphi<0$ in $\Omega$ and so {\bf (WMP)} fails in $\Omega$.

For the remaining case $\lambda \geq 0$ and $\mu < 0$, there exists $(\lambda_1, \mu_1)\in \mathcal{C}_1$ with $\lambda_1 > 0$ large enough (and so $\mu_1>0$ small enough) so that $\lambda > -\lambda_1$ and $\mu < -\mu_1$. Therefore,

\[
\left\{
\begin{array}{llll}
-\Delta_p \varphi - \lambda a(x) \vert -\psi\vert^{\beta_1-1}(-\psi) &=& \lambda_1 a(x) \psi^{\beta_1} + \lambda a(x)\psi^{\beta_1} \\
&=& (\lambda + \lambda_1)a(x) \psi^{\beta_1} \geq (\not\equiv)\ 0\ \ & {\rm in}\ \Omega; \\
-\Delta_q (-\psi) - \mu b(x) \varphi^{\beta_2} &=& - \mu_1 b(x) \varphi^{\beta_2} - \mu b(x) \varphi^{\beta_2} \\
&=& -(\mu + \mu_1)b(x) \varphi^{\beta_2} \geq (\not\equiv)\ 0\ \ & {\rm in}\ \Omega
\end{array}
\right. 
\]
and $\varphi = 0 =-\psi$ on $\partial\Omega$. But, $-\psi<0$ in $\Omega$ and so again {\bf (WMP)} fails in $\Omega$.

Conversely, we next show that {\bf (WMP)} holds in $\Omega$ for any couple $(\lambda, \mu) \in \overline{\mathcal{R}_1} \setminus \mathcal{C}_1$. Since $\Delta_p$ and $\Delta_q$ satisfies weak maximum principle in $\Omega$, we have {\bf (WMP)} holds in $\Omega$ if either $\lambda = 0$ and $\mu \geq 0$ or $\lambda \geq 0$ and $\mu = 0$. Now let $(\lambda,\mu) \in \mathcal{R}_1$. Let $(u,v)$ be a weak solution of the system (\ref{supersol}).  Note that $\lambda < \lambda_1$ and $\mu < \mu_1$, where $(\lambda_1, \mu_1)$ is a principal eigenvalue of (\ref{1.3}) with $\frac{\mu}{\lambda}=\frac{\mu_1}{\lambda_1}$. Let $(\varphi,\psi)$ be a positive eigenfunction associated to $(\lambda_1, \mu_1)$. Arguing by contradiction, assume that $u$ or $v$ is negative somewhere in $\Omega$. Then, by strong comparison principle and Hopf's lemma for the $p$-Laplacian (see \cite{GS, 23, 25}), there exists some $\gamma>0$ such that 
\[
-u\leq\gamma\varphi\ \text{and}\ -v\leq\gamma^\omega\psi\ \ \text{in}\ \ \Omega,
\]
where $\omega:=\frac{p-1}{\beta_1}$. Let $\overline{\gamma}$ be the minimum of such $\gamma's$. Thus, $\overline{\gamma}>0$. Since $\lambda < \lambda_1$ and $\mu < \mu_1$, we derive
\[
\left\{
\begin{array}{llll}
-\Delta_p (-u) \leq \lambda a(x) \vert -v\vert^{\beta_1-1}(-v)\leq  \lambda a(x) (\overline{\gamma}^\omega\psi)^{\beta_1}\leq (\not\equiv)\ \lambda_1 a(x)  (\overline{\gamma}^\omega\psi)^{\beta_1}=-\Delta_p(\overline{\gamma}\varphi)\ \ & {\rm in}\ \Omega; \\
-\Delta_q (-v) \leq \mu b(x) \vert -u\vert^{\beta_2-1}(-u)\leq  \mu b(x) (\overline{\gamma}\varphi)^{\beta_2}\leq (\not\equiv)\ \mu_1 b(x)  (\overline{\gamma}\varphi)^{\beta_2}=-\Delta_q(\overline{\gamma}^\omega\psi)\ \ & {\rm in}\ \Omega
\end{array}
\right.
\]
and $\overline{\gamma}\varphi= -u =\overline{\gamma}^\omega\psi=- v = 0$ on $\partial \Omega$. It follows from the strong comparison principle to each above equation (see Theorem A.1 of \cite{Cuesta}) that $ -u< \overline{\gamma}\varphi$ and $-v<\overline{\gamma}^\omega\psi$ in $\Omega$. Thus, we can find $0<\varepsilon <1$ such that $-u\leq\varepsilon\overline{\gamma}\varphi$ and $-v\leq(\varepsilon\overline{\gamma})^\omega\psi$ in $\Omega$, contradicting the definition of $\overline{\gamma}$. Then, $u, v \geq 0$ in $\Omega$. Hence, we complete the wished proof of theorem. \ \rule {1.5mm}{1.5mm}

\section{Proof of Theorem \ref{CP}}

In this section we establish the characterization of {\bf (WCP)} and {\bf (SCP)} in terms of the principal curve $\mathcal{C}_1$ as stated in Theorem \ref{CP}. Using strong comparison principle (see Theorem A.1 of \cite{Cuesta}) and arguing in a similar way as in the proof of Theorem \ref{MP}, we see that it suffices to show only that (i) $\Leftrightarrow$ (ii).

(ii) $\Rightarrow$ (i) Taking $u, v \equiv 0$ in $\Omega$ (and so $f_1,f_2\equiv 0$ in $\Omega$), we have {\bf (WCP)} in $\Omega$ implies {\bf (WMP)} in $\Omega$. Then, by Theorem \ref{MP}, {\bf (WCP)} in $\Omega$ implies $(\lambda, \mu) \in \overline{\mathcal{R}_1} \setminus \mathcal{C}_1$. 

(i) $\Rightarrow$ (ii) Consider $(\lambda, \mu) \in \overline{\mathcal{R}_1} \setminus \mathcal{C}_1$. Since $\Delta_p$ satisfies weak maximum and comparison principles in $\Omega$, the conclusion is direct in the cases that $\lambda = 0$ or $\mu = 0$. Thus, it suffices to consider $(\lambda, \mu) \in \mathcal{R}_1$. In this case, since $(u,v)$ is a weak solution of the problem \eqref{comp1}, by Theorem \ref{MP}, if $f_1,g_1\equiv 0$ in $\Omega$ we have $u, v \equiv 0$ in $\Omega$ and if $f_1+g_1\not\equiv 0$ in $\Omega$, we get $u, v > 0$ in $\Omega$. Note that, if $f_1, g_1 \equiv 0$ in $\Omega$, then the conclusion follows readily from {\bf (WMP)} in $\Omega$. Assume then $f_1+g_1\not\equiv 0$ in $\Omega$ (and so $u, v > 0$ in $\Omega$). Thus, we also have $f_2+g_2\not\equiv 0$ in $\Omega$. Therefore, since $(z,w)$ is a weak solution of \eqref{comp2}, by {\bf (SMP)}, we have $z, w > 0$ in $\Omega$. 

To finish the proof of {\bf (WCP)}, it suffices to show that $u \leq z$ in $\Omega$. This follows directly by using the fact that $\Delta_q$ satisfies weak comparison principle in $\Omega$ (see \cite{23}). Assume by contradiction that $u > z$ somewhere in $\Omega$. In this case, the set $\Gamma := \{\gamma>0 : z > \gamma u\ {\rm and} \ w > \gamma^{\omega} v\ {\rm in} \ \Omega\}$, where $\omega:=\frac{p-1}{\beta_1}$, is nonempty by Hopf's Lemma (see \cite{25}) and is also upper bounded. Set $\overline{\gamma}:=\sup \Gamma > 0$. Notice that $z \geq \overline{\gamma} u$ and $w \geq \overline{\gamma}^{\omega} v$ in $\Omega$. Note also that the statement of contradiction leads to $\overline{\gamma} < 1$. Thus, since $f_1+g_1\not\equiv 0$ and $f_2+g_2\not\equiv 0$ in $\Omega$, we derive

\[
\left\{
\begin{array}{lll}
-\Delta_p (\overline{\gamma}u) =  \lambda a(x) (\overline{\gamma}^\omega v)^{\beta_1}+\overline{\gamma}^{p-1}f_1\leq (\not\equiv)\ \lambda a(x)  w^{\beta_1}+f_2=-\Delta_p(z) \ \ & {\rm in}\ \Omega; \\
-\Delta_q (\overline{\gamma}^\omega v) =  \mu b(x) (\overline{\gamma}u)^{\beta_2}+\overline{\gamma}^{\beta_2}g_1\leq (\not\equiv)\ \mu b(x)  z^{\beta_2}+g_2=-\Delta_q(w) \ \ & {\rm in}\ \Omega
\end{array}
\right. 
\]
and $\overline{\gamma}u= z =\overline{\gamma}^\omega v=w = 0$ on $\partial \Omega$. It follows from the strong comparison principle to each above equation (see Theorem A.1 of \cite{Cuesta}) that $ z> \overline{\gamma} u$ and $w>\overline{\gamma}^\omega v$ in $\Omega$. Thus, we can find $0<\varepsilon <1$ such that $z \geq (\overline{\gamma} + \varepsilon) u$ and $w \geq (\overline{\gamma}+ \varepsilon)^\omega v$ in $\Omega$, contradicting the definition of $\overline{\gamma}$. This concludes the desired proof. \ \rule {1.5mm}{1.5mm}

\section{Proof of Theorem \ref{lower}} 

Let $a,b\in L^\infty(\Omega) \cap C(\Omega)$, $(\lambda_1, \mu_1)\in \mathcal{C}_1$ and $(\varphi,\psi) \in X$ be a positive eigenfunction of the problem \eqref{1.3} corresponding to $(\lambda_1, \mu_1)$. Since $\lambda_1 a\varphi^{\beta_1},\mu_1 b\psi^{\beta_2}\in C(\Omega)$, by Theorem 1.8 of \cite{medina}, we have $(\varphi,\psi)$ is a viscosity subsolution of the system (\ref{1.3}). Thus, applying the ABP estimate for the $p$-Laplacian (see Theorem 3 of \cite{ACI}) to the first equation of \eqref{1.3}, we have
\begin{eqnarray}\label{i1}
||\varphi||_{L^\infty(\Omega)} &=& \sup_\Omega \varphi \leq c(n,p)d\lambda_1^{\frac{1}{p-1}} \|a\|^{\frac{1}{p-1}}_{L^\infty(\Omega)} \|\psi\|^{\frac{\beta_1}{p-1}}_{L^\infty(\Omega)}\vert\Omega\vert^{\frac{1}{n(p-1)}}
\end{eqnarray}
and for the $q$-Laplacian to the second equation of \eqref{1.3}, we obtain
\begin{eqnarray}\label{i2}
||\psi||_{L^\infty(\Omega)} \leq c(n,q)d\mu_1^{\frac{1}{q-1}} \|b\|^{\frac{1}{q-1}}_{L^\infty(\Omega)} \|\varphi\|^{\frac{\beta_2}{q-1}}_{L^\infty(\Omega)}\vert\Omega\vert^{\frac{1}{n(q-1)}}.
\end{eqnarray}

Therefore, joining inequalities (\ref{i1}) and (\ref{i2}) and using that 
\[
\beta_1 \beta_2 = (p-1)(q-1)\ \text{ and }\ \lambda_1^{\frac{1}{r}}\mu_1^{\frac{1}{s}} = \Lambda',
\]
we finally derive \eqref{lb1} and conclude the proof. \ \rule {1.5mm}{1.5mm}

\section{Proof of Theorem \ref{sm}}

By Theorem \ref{MP}, the necessity of (i) and the equivalence between {\bf (WMP)} and {\bf (SMP)} follow directly. Thus, it suffices to prove that the assertion (i) implies (ii).

Assume that $\lambda \geq 0$ and $\mu \geq 0$. If either $\lambda = 0$ or $\mu = 0$, then by Theorem \ref{MP}, the desired {\bf (WMP)} follow.

Finally, assume now that $\lambda > 0$ and $\mu > 0$. We consider the positive constant $\eta$ given by

\[
\eta := \frac{1}{\left[\lambda^{\frac{1}{r}}\mu^{\frac{1}{s}}\left(c(n,p)^{\frac{\beta_2}{s}}c(n,q)^{\frac{\beta_1}{r}}d^{\frac{\beta_2}{s}+\frac{\beta_1}{r}} \|a\|^{\frac{1}{r}}_{L^\infty(\Omega)} \|b\|^{\frac{1}{s}}_{L^\infty(\Omega)}\right)\right]^{\frac{nrs}{r+s}}},
\]
where $d:=\operatorname{diam }(\Omega)$, $r:=\sqrt{\beta_1(p-1)}$, $s:=\sqrt{\beta_2(q-1)}$ and $c(n,p)$ and $c(n,q)$ are the explicit constants of ABP estimate for the $p$-Laplacian and $q$-Laplacian, respectively. Then, by using the estimate \eqref{lb1} of Theorem \ref{lower}, we obtain

\[
\Lambda' \geq \frac{\vert\Omega\vert^{-\left(\frac{1}{ns}+\frac{1}{nr}\right)}}{c(n,p)^{\frac{\beta_2}{s}}c(n,q)^{\frac{\beta_1}{r}}d^{\frac{\beta_2}{s}+\frac{\beta_1}{r}} \|a\|^{\frac{1}{r}}_{L^\infty(\Omega)} \|b\|^{\frac{1}{s}}_{L^\infty(\Omega)}}>\lambda^{\frac{1}{r}}\mu^{\frac{1}{s}}
\]
whenever $|\Omega| < \eta$. Therefore, we get $(\lambda, \mu) \in \mathcal{R}_1$ for such domains and so, by Theorem \ref{MP} the assertion (ii) holds. This concludes the proof of theorem.\ \rule {1.5mm}{1.5mm}

\section{Proof of Theorem \ref{FHA}}

Let $(\lambda,\mu) \in \overline{\mathcal{R}_1} \setminus \mathcal{C}_1$ and $f,g\in L^{\infty}(\Omega)$ such that $f,g\geq 0$ in $\Omega$. If either $\lambda = 0$ or $\mu = 0$, then clearly, by existence and uniqueness of weak solution for the problem (\ref{P1}), the system (\ref{supersol}) admits a unique weak solution $(u,v) \in W^{1,p}_0(\Omega)\times W^{1,q}_0(\Omega)$. Applying the weak maximum principle to each equation of (\ref{supersol}), we get $u, v \geq 0$ in $\Omega$.

Assume now that $(\lambda,\mu) \in \mathcal{R}_1$. Then, by Theorem 3.1 of \cite{Cuesta}, the system (\ref{supersol}) admits a unique weak solution $(u,v) \in X_+$. Let $(z,w)\in W^{1,p}_0(\Omega)\times W^{1,q}_0(\Omega)$ be a weak solution of the system (\ref{supersol}). So, by {\bf (WCP)}, $u=z$ and $v=w$ in $\Omega$. Therefore, the system (\ref{supersol}) admits a unique weak solution in $W^{1,p}_0(\Omega)\times W^{1,q}_0(\Omega)$ and this is nonnegative.

Conversely, assume that $(\lambda,\mu) \in \R^2$ and the problem \eqref{supersol} admits a unique weak solution $(u,v) \in W^{1,p}_0(\Omega)\times W^{1,q}_0(\Omega)$ and that this satisfies $u, v \geq 0$ in $\Omega$ for any pair $(f,g) \in (L^\infty(\Omega))^2$ of nonnegative functions. Thus, {\bf (WMP)} associated to \eqref{1.3} holds in $\Omega$ and so by Theorem \ref{MP}, we have $(\lambda,\mu) \in \overline{\mathcal{R}_1} \setminus \mathcal{C}_1$. This ends the proof. \ \rule {1.5mm}{1.5mm}\\

%\n {\bf Acknowledgments:} The author was partially supported by Fapemig (Universal-APQ-00709-18). The author is indebted to the anonymous referee for his/her valuable comments.


\begin{thebibliography}{99}

\addcontentsline{toc}{section}{References}

\bibitem{AFM}
{B. Alziary, J. Fleckinger, M.-H. Lécureux} - \textit{Principal eigenvalue and Maximum principle for some elliptic systems defined on general domains with refined Dirichlet boundary condition}, Comm. Math. Anal. 7 (2009), 1-11.

\bibitem{Am2}
{H. Amann} - \textit{Maximum principles and principal eigenvalues, in 10 Mathematical Essays on Approximation in Analysis and Topology}, (J. Ferrera, J. López-Gómez and F. R. Ruiz del Portal eds.), pp. 1-60, Elsevier, Amsterdam, 2005.

\bibitem{AnLo}
{I. Anton, J. López-Gómez} - \textit{Principal eigenvalue and maximum principle for cooperative periodic-parabolic systems}, Nonlinear Analysis 178 (2019), 152-189.

\bibitem{ACI}
{R. Argiolas, F. Charro, I. Peral} - \textit{On the Aleksandrov-Bakel'man-Pucci Estimate for Some Elliptic and Parabolic Nonlinear Operators}, Arch. Rational Mech. Anal. 202 (2011), 875-917.

\bibitem{BeNiVa}
{H. Berestycki, L. Nirenberg, S.R.S. Varadhan} - \textit{The principal eigenvalue and maximum principle for second order elliptic operators in general domains}, Comm. pure Appl. Math. (1994), 47-92.

\bibitem{CaLo}
{S. Cano-Casanova, J. López-Gómez} - \textit{Properties of the principal eigenvalues of a general class of nonclassical mixed boundary value problems}, J. Diff. Eq. 178 (2002), 123-211.

\bibitem{CauLo}
{P.A. Caudevilla, J. López-Gómez} - \textit{Asymptotic behaviour of principal eigenvalues for a class of cooperative systems}, J. Diff. Eq. 244 (2008), 1093-1113.

\bibitem{ClDfMi}
{Ph. Cl\'{e}ment, D. G. de Figueiredo, E. Mitidieri} - \textit{Positive solutions of semilinear elliptic systems}, Comm. in PDE 17 (1992), 923-940.

\bibitem{2}
{Ph. Cl\'{e}ment, R.F. Manásevich, E. Mitidieri} - \textit{Positive solutions for a quasilinear
system via blow up}, Comm. P.D.E., 18 (1993), 2071-2106.

\bibitem{CT}
{M. Cuesta, P. Tak\'{a}c} - \textit{A strong comparison principle for the Dirichlet $p$-Laplacian}, in
Reaction Diffusion Systems, (G. Caristi and E. Mitidieri, Eds.). Lecture Notes in Pure
and Applied Mathematics, 194. Marcel Dekker (1997).

\bibitem{Cuesta}
{M. Cuesta, P. Tak\'{a}c} - \textit{Nonlinear eigenvalue problems for degenerate elliptic systems}, Differ. Integral Equations 23 (2010), 1117-1138. 

\bibitem{Da}
{L. Damascelli} - \textit{Comparison theorems for some quasilinear degenerate elliptic operators
and applications to symmetry and monotonicity results}, Ann. Inst. H. Poincaré. Analyse non linéaire 15 (1998), 493-516.

\bibitem{DaSc}
{L. Damascelli, B. Sciunzi} - \textit{Regularity, monotonicity and symmetry of positive
solutions of $m$-Laplace equations}, J. Differential Equations 206 (2004), 483-515.

\bibitem{DaSc1}
{L. Damascelli, B. Sciunzi} - \textit{Harnack inequalities, maximum and comparison
principles, and regularity of positive solutions of $m$-Laplace equations}, Calc. Var.
Partial Differential Equations 25 (2006), 139-159.

\bibitem{DfFe}
{D. G. de Figueiredo, P. L. Felmer} - \textit{On superquadratic elliptic systems}, Trans. Amer. Math. Soc. 343 (1994), 99-116.

\bibitem{7}
{DE. DiBenedetto} - \textit{$C^{1+\alpha}$ local regularity of weak solutions of degenerate elliptic equations}, Nonlinear Anal. 7 (1983), 827-850.

\bibitem{faber}
{G. Faber} - \textit{Beweis, dass unter allen homogenen Membranen von gleicher Fl\"{a}che und
gleicher Spannung die kreisf\"{o}rmige den tiefsten Grundton gibt (in German)}, Sitzungberichte
der mathematisch-physikalischen Klasse der Bayerischen Akademie der Wissenschaften
zu München Jahrgang (1923), 169-172.

\bibitem{FeMa}
{P. Felmer, S. Mart\'{i}nez} - \textit{Existence and uniqueness of positive solutions to certain differential systems},
Adv. Differential Equations 4 (1998), 575-593.

\bibitem{FGT}
{J. Fleckinger-Pellé, J. P. Gossez, P. Takác, F. de Thélin} - \textit{Nonexistence of solutions and an anti-maximum principle for cooperative systems with the $p$-Laplacian}, Math. Nachr. 194 (1998), 49-78.

\bibitem{FHTT}
{J. Fleckinger-Pellé, J. Hernández, P. Takác, F. de Thélin} - \textit{Uniqueness and positivity for solutions of equations with the $p$-Laplacian}. Reaction diffusion systems (Trieste, 1995), 141-155, Lecture Notes in Pure and Appl. Math., 194, Dekker, New York, 1998.

\bibitem{FT}
{J. Fleckinger-Pellé, P. Takác} - \textit{Uniqueness of positive solutions for nonlinear cooperative systems with the $p$-Laplacian}, Indiana Univ. Math. J. 43 (1994), 1227-1253.

\bibitem{GS}
{J. García-Melián, J. C. Sabina de Lis} - \textit{Maximum and comparison principles for operators involving the $p$-Laplacian}, J. Math. Anal. Appl. 218 (1998), 49-65.

\bibitem{GV}
{M. Guedda, L. Veron} - \textit{Quasilinear elliptic equations involving critical Sobolev
exponents}, Nonlinear Anal. 13 (1989), 879-902.

\bibitem{HuvdV}
{J. Hulshof, R. C. A. M. van der Vorst} - \textit{Differential systems with strongly indefinite variational structure}, J. Funct. Analysis, 114 (1993), 32-58.

\bibitem{krahn}
{E. Krahn} - \textit{\"{U}ber eine von Rayleigh formulierte Minimaleigenschaft des Kreises}, Math. Ann. 91 (1925), 97-100.

\bibitem{LM}
{E. J. F. Leite, M. Montenegro} - \textit{Maximum and comparison principles to Lane-Emden systems}, Journal of the London Mathematical Society 101 (2020), 23-42.

\bibitem{EM1}
{E. J. F. Leite, M. Montenegro} - \textit{Principal curves to nonlocal Lane-Emden systems and related
	maximum principles}, Calc. Var. Partial Differential Equations 59 (2020), 118.


\bibitem{L}
{E. J. F. Leite} - \textit{On the principal eigenvalues of the degenerate elliptic systems}, Electronic Journal of Qualitative Theory of Differential Equations 40 (2020), 1-15.

\bibitem{15}
{G. Lieberman} - \textit{Boundary regularity for solutions of degenerate elliptic equations}, Non-linear Anal. 12 (1988), 1203-1219.

\bibitem{Lo1}
{J. López-Gómez} - \textit{Linear Second Order Elliptic Operators}, World Scientific, Singapore, 2013.

\bibitem{Lo}
{J. López-Gómez} - \textit{The maximum principle and the existence of principal eigenvalue for some linear weighted boundary value problems}, J. Diff. Eq. 127 (1996), 263-294.
%
\bibitem{LoMo}
{J. López-Gómez, M. Molina-Meyer} - \textit{The maximum principle for cooperative weakly coupled elliptic systems and some applications}, Diff. Int. Equations 7 (1994), 383-398.

\bibitem{medina}
{M. Medina, P. Ochoa} - \textit{On viscosity and weak solutions for non-homogeneous $p$-Laplace equations}, Adv. Nonlinear Anal. 8 (2019), 468-481.

\bibitem{Mi1}
{E. Mitidieri} - \textit{A Rellich type identity and applications}, Comm. Partial Differential Equations 18 (1993), 125-151.

\bibitem{Mo4}
{M. Montenegro} - \textit{The construction of principal spectra curves for Lane-Emden systems and applications}, Ann. Scuola Norm. Sup. Pisa Cl. Sci. 29 (2000), 193-229.

\bibitem{PS}
{P. Pucci, J. Serrin} - \textit{The maximum principle}, Birkh\"{a}user Verlag, Basel, 2007.

\bibitem{SZ1}
{J. Serrin, H. Zou} - \textit{Existence of positive entire solutions of elliptic Hamiltonian systems}, Comm. Partial Differential Equations 23 (1998), 577-599.


\bibitem{Sw}
{G. Sweers} - \textit{Strong positivity in $C(\overline\Omega)$ for elliptic systems}, Math. Z. 209 (1992), 251-271.

\bibitem{23}
{P. Tolksdorf} - \textit{On the Dirichlet problem for quasilinear equations in domains with conical boundary points}, Comm. P.D.E., 8 (1983), 773-817.

\bibitem{24}
{P. Tolksdorf} - \textit{Regularity for a more general class of quasilinear elliptic equations}, J.
Differential Equations, 51 (1984), 126-150.

\bibitem{25}
{J. L. Vázquez} - \textit{A strong maximum principle for some quasilinear elliptic equations}, Appl. Math. Optim., 12 (1984), 191-202.

\end{thebibliography}
\end{document}